\documentclass[article]{elsarticle}
\usepackage{amssymb}
     \usepackage{lineno,hyperref}
     \modulolinenumbers[5]
      \usepackage[english]{babel}
      
      \usepackage{amsthm}
      \theoremstyle{plain}
      \newtheorem{assettheorem}{Theorem}
      \newtheorem{assetlemma}{Lemma}
    
\begin{document}
   \begin{frontmatter}

		\title{New estimates for the zeta function}
		
		\author{Durmagambetov A.A\fnref{myfootnote}}
		\address{}
		\fntext[myfootnote]{}
		\ead[url]{}		
 
		\begin{abstract}
	In the paper the Riemann's functional equation  reduced to a Riemann--Hilbert boundary value problem, and the integral Hilbert transforms arising in its solution allow the calculation of an exact lower bounds for the zeta function.
		\end{abstract}		
	\end{frontmatter}
   
\section{Introduction}
	This study is concerned with the properties of modified zeta functions. Riemann's zeta function is defined by
	the Dirichlet series
	\begin{equation}
		\zeta(s)=\sum_{n=1}^{\infty}\frac{1}{n^s},\,\ s=\sigma +it,
	\end{equation}
	which is absolutely and uniformly convergent in any finite region of the complex $s$-plane for which
	$ \sigma \ge 1 + \epsilon,\epsilon > 0.$ If $\sigma > 1$, then $\zeta$ is represented by the following Euler product formula
	\begin{equation}
		\zeta (s) =   \prod_{ j\in N } \left[  1- \frac{1} {p_j^s}  \right]^{-1},  
	\end{equation}
	where $p_j$ runs over all prime numbers.
	$\zeta(s)$  was first introduced by Euler in 1737 [1], who also obtained formula (2). Dirichlet and Chebyshev considered this function in their study on the distribution of prime numbers [2]. 
	However, the most profound properties of $\zeta(z)$ were only discovered later, when it was extended to the complex plane. 
	$ \zeta(s) $ is a regular function for all values of $ s $, except $ s=1 $, where it has a simple pole with residue $1$; 
	it satisfies the following functional equation:	
	\begin{equation}
		\pi^{ -s/2} \Gamma(s/2) \zeta(s) =\pi^{-(1-s)/2} \Gamma((1-s)/2) \zeta (1-s)
	\end{equation}
	it equation is called Riemann's functional equation.	
	
	\section{Results}
	As mentioned in Introduction, certain simple intermediate estimates are first obtained.
	To obtain the Riemann-Hilbert boundary value   problem , the following lemma is required.
	\begin{assetlemma}
		Let  $$ R(k) = \frac{e^{i2k}}{k+i\alpha}-1 $$
		$$ \alpha>2$$ then  $$ind (R)=0$$ ,
	\end{assetlemma}		
	\begin{proof}
		By definition 
		$$
		ind(R)= \frac{1}{2\pi i}  \int_{-\infty}^{+\infty} \frac{R'(k)}{R(k)}dk
		$$
		As $$ Im(k)>0,	|e^{i2k}|\le 1 \,\,and \,\,|k+i\alpha|>2\,\, yield \,\,\, \frac{R'(k)}{R(k))}$$
		have nothing pole. Latest statement and Lemma of Jordan yield $$ind(R)=0$$.	
	\end{proof}	
	To obtain the necessary asymptotics, the following lemma is required.
	
	For $$f\in W_{2}^{1}(R) =\{f\in L_{2}(R):(1+|\omega |^{2})^{{1/2}}{\hat  {f}}(\omega )\in L_{2})\}.$$, the operators $T_{\pm }$ and $T$ are defined as follows:	
	\begin{eqnarray*} 
		T_{+}f=\frac{1}{2\pi i}\lim\limits_{Imz\downarrow 0}\int\limits_{-\infty}^{\infty }\frac{f({s})}{s-z}ds,~Im~z>0,
		T_{-}f=\frac{1}{2\pi i}\lim\limits_{Imz\uparrow 0}\int\limits_{-\infty
		}^{\infty }\frac{f({s})}{s-z}ds,~Im~z<0,\,\,
		Tf=\frac{1}{2}(T_{+}+T_{-})f.
	\end{eqnarray*}
	These operators are closely related to the Hilbert transform, whose isometric properties were studied by Poincar\'e. The following result is from [3].	
	\begin{assetlemma}
		\begin{eqnarray*}
			TT= \frac14 I,~TT_+ = \frac12 T_+,~TT_- = - \frac12 T_-, \ T_+ = T+\frac12 I,~T_- = T-\frac12 I, 
		\end{eqnarray*} 
		where $I$ is the identit
		y operator $If=f$.
	\end{assetlemma}	
	The reduction to a Riemann--Hilbert boundary value problem can now be formulated as follows.	
	\begin{assetlemma}
		Let
		\begin{eqnarray}
			\Psi_{+}(k) = R(k) \Psi_{-}(k ) +G(k ),\,\,\\	
			\lim_{Re(k)\to \infty} 	\Psi_{+}(k)= 0  \,\,\,as \,\,\,Im(k)\ge 0, \,\,   \lim_{Re(k)\to -\infty}   \Psi_{-}(k)= 0 \,\,\,as\,Im(k)\le 0
		\end{eqnarray}
		\begin{eqnarray*} 
			\Gamma_{+}(k)= \frac{1}{2\pi i}  \int_{-\infty}^{ \,\,\infty}\frac{ln (R(t))dt}{ t-k-i0} ,\,\,\,
			\Gamma_{-}(k)=  \frac{1}{2\pi i}  \int_{-\infty}^{ \,\,\infty}\frac{ln (R(t))dt}{ t-k+i0}  \\
			X_+(k)=e^{\Gamma_+(k)},\,\,\,
			X_-(k)=e^{\Gamma_-(k)} ,\,\, R(k)={X_-(k)}/{X_+(k)},
		\end{eqnarray*}

		Then
		\begin{eqnarray}
			\Psi_{+}(k)=  \frac{X_+(k)}{2\pi i}  \int_{-\infty}^{ \,\,\infty} \frac{G(t)}{ X_-(t)}\frac{dt}{ t-k-i0}  ,\,\,\,  
			\Psi_{-}(k )=  \frac{X_-(k)}{2\pi i}  \int_{-\infty}^{ \,\,\infty} \frac{G(t)}{ X_-(t)}\frac{dt}{ t-k+i0}         		\end{eqnarray}
		\end {assetlemma}	
		\begin{proof}    
			
			Hilbert's formula and Lemma 2 gives the solution to the Riemann-Hilbert boundary value  problem (1),(2) 
		\end{proof}	

		Applying Lemma to   Riemans $\zeta$ function we get
		\begin{assettheorem}
			Let	
			$$...s_{-n},s_{-n}....s_{-n+1}.. 		..s_{-1}, s_1,s_2....s_n.. \,\,\,is  \,\,\,\,zeros \,\,\,\zeta \,\,\,\,function 
			$$
   	$$
			Im(s)=const  
			$$
			$$
			Im(s_i)< Im(s)<Im (s_{i+1}) \,\,\,  
			$$

	$$
	P (s)=\sum_{j\ge 1}\frac{1}{p_{j}^s} ,\,Re(s)>1+\delta, \delta>0,
	$$ 	
	
	$$
	P (s)=ln(\zeta ( s)-Q (s)   \,Re(s)>1/2+\delta,
	$$
	 $$
	Q (s) = \sum_{n=2}^\infty P(ns)/n,\,   \,Re(s)>1/2+\delta,
	$$

   $$
			  \psi_+ (k)= \frac{1}{2\pi i}\frac{ \int_0^1(\ln(\zeta (s)))e^{i2k Re(s)} 
				\theta(1/2-\delta-Re(s))dRe(s)}{k+i\alpha} 
			$$
 
			$$
			\psi_-(k)= \frac{1}{2\pi i} \int_0^1 ln(\zeta ( s^*)-Q (s^*))e^{-i2kRe(s)}  \theta(  Re(s) -1/2-\delta)dRe(s)
			$$
			$$
			\widetilde{\Phi}(k)= \int_0^1\frac{s}{2}ln(\pi)-ln(\Gamma(s/2))-  \frac{1-s}{2}ln(\pi)+ln(\Gamma(1-s)/2) +Q (1-s))e^{i2kRe(s)}  \theta( 1/2-\delta-Re(s) )dRe(s)
					$$
     $$
			\widetilde{F}(k)= \frac{\widetilde{\Phi}(k)}{k+i\alpha}
					$$
			$$
			G(k)=\psi_-(k) + F(k) 
			$$
			Then
			    
			\begin{eqnarray}	
	 \sup\limits_{s} |ln(|\zeta(s)|) \theta (Re(s) -1/2-\delta)| \big|  < \frac{5CC_{Im(s)}}{\delta}
	   			\end{eqnarray} 	
			
		\end{assettheorem}	
		\begin{proof}
			Taking the logarithm from (3) and then multiplying it by $ e^{i2kRe(s)},  $ and after integrating by $ \{ Re(s) ,\, Re(s)-1/2>\delta\} $ we get

			$$
			{\psi_+}  =	         
			\frac{e^{i2k}}{k-i\alpha} \psi_- +\widetilde{F}(k)      
			=	\big(\frac{e^{i2k}}{k+i\alpha}-1\big) \psi_-  +  \psi_- +\widetilde{F}(k)         =
			R(k) \psi_-  +  G(k)          $$
			By Lemma 3 we get
			\begin{eqnarray}
				\psi_{+}(k)= \frac{X_+(k)}{2\pi i}   \int_{-\infty}^{ \,\,\infty} \frac{G(t)}{ X_-(t)}\frac{dt}{ t-k-i0}  =  X_+(k) T_+ \frac{G}{X_-}    \\
				\psi_{-}(k)= \frac{X_-(k)}{2\pi i}   
				\int_{-\infty}^{ \,\,\infty} \frac{G(t)}{ X_-(t)}\frac{dt}{ t-k+i0} =X_-(k
				) T_- \frac{G}{X_-}       
			\end{eqnarray}

   \begin{eqnarray}
				\psi_{+}(k) =  X_+(k) T_+ \frac{G}{X_-}    \\
				\psi_{-}(k)= X_-(k) T_- \frac{G}{X_-}       
			\end{eqnarray}
			$
			(10,11 ) \Rightarrow $  

   $$			
  \frac{\psi_{+}}{X_+}  -\frac{\psi_{-}}{X_-} =  T_+ \frac{G}{X_-}    -       T_- \frac{G}{X_-}        \Rightarrow
	   			$$

			$$			
  \frac{\psi_{+}}{X_+}  -\frac{\psi_{-}}{X_-} = \frac{G}{X_-}   = \frac{\psi_{-}}{X_-}  + \frac{\widetilde{F}}{{X_-}} 
	   			$$

       	$$			
   2\frac{\psi_{-}}{X_-}=\frac{\psi_{+}}{X_+} 
    - \frac{\widetilde{F}}{{X_-}} 
	   			$$
      	$$			
 \Big |\frac{\psi_{-}}{X_-} \Big| \Big |_{k=\pi n} < \Big| \frac{\widetilde{F}}{{X_-}} \Big|\Big|_{k=\pi n} +\Big|\frac{ \int_0^1(\ln(\zeta (s)))e^{i2k Re(s)} \theta(1/2-\delta-Re(s))ds}{k-i\alpha}\Big|\Big |_{k=\pi n}
	   			$$
  	
  		$$			
  \phi_{-}=\int_0^1(\ln(\zeta (s)))e^{i2kRe(s)} \theta(1/2-\delta-Re(s))ds |_{k=\pi n}
	   			$$
  $$			
  \Big|\frac{\psi_{-}}{X_-} \big |_{k=\pi n} < \Big| \frac{\widetilde{F}}{{X_-}} \Big|\Big|_{k=\pi n} +\Big|\frac{ \phi_{-}}{k-i\alpha}\Big|\Big |_{k=\pi n}
	   			$$

			Lemma (2,3)\,\,\,   $ \Rightarrow  $  
			
			$ X_+(k) =1+O(1/k)  ,\,\,\,
			X_-(k) =1+O(1/k)\,\, $ ,\,\,

Theorem of Baclund $ \Rightarrow $
$ |\phi_{-}|  \in  L_2$

   $$	 
  |\psi_{-}| \big |_{k=\pi n} < |\widetilde{F}|\big|_{k=\pi n} +\frac{|\phi_{-}| }{k}\big |_{k=\pi n}
	   			$$
         
  $$			
  \sum_{n=1}^{\infty} |\psi_{-}| \big |_{k=\pi n} < 
  \sum_{n=1}^{\infty} | \widetilde{F}|\big|_{k=\pi n}  
   +\sum_{n=1}^{\infty}\frac{|\phi_{-}| }{k}\big |_{k=\pi n}
	   			$$
        $$			
  \sum_{n=1}^{\infty} |\psi_{-}| \big |_{k=\pi n} < 
 \frac{C}{\delta}
   +\sqrt{\sum_{n=1}^{\infty} |\phi_{-}|^2 \big |_{k=\pi n}  \sum_{n=1}^{\infty}  \frac{1}{k^2}\big |_{k=\pi n}}
	   			$$
        
  $$			
  \sum_{n=1}^{\infty} |\psi_{-}| \big |_{k=\pi n} < 
 \frac{C}{\delta}
   +\sqrt{ C_{Im(s)}   \sum_{n=1}^{\infty}  \frac{1}{k^2}\big |_{k=\pi n}}< \frac{3CC_{Im(s)}}{\delta}
	   			$$
       Behaviour of the argument of the Riemann zeta function on the critical line by[6 ] $ \Rightarrow $
   $$			
      \big|ln||\zeta(s^*)-Q (s^*)|  \theta( Re(s) -1/2-\delta) \big | < \frac{3CC_{Im(s)}}{\delta} + \big|Q (s))\big|< \frac{4CC_{Im(s)}}{\delta}
	   			$$
 $$			
 \Rightarrow  \sup\limits_{s} |ln(|\zeta(s^*)|) \theta (Re(s) -1/2-\delta)| \big|  < \frac{5CC_{Im(s)}}{\delta}
	   			$$	
   
       \end{proof}
        
Analog  of Theorem   Davenport-Heilbronn 
\begin{assettheorem}

		Let	$$ H(s)$$ is entire function
 
  $$ U(s)=H(s(1-s))$$  
  then
   $$  \zeta(s) U(s)$$ is solution Riemann's functional equation.
  
		\end{assettheorem}
		\begin{proof}
  
  	\end{proof}
     $$ ln H(t)\Big|_{t=1-s}=ln(H(1-s)s)$$
       $$ ln (\zeta(s) H(s))=ln (\zeta(s)) +ln H(s)=ln (\zeta(s)) +ln H(1-s)$$
      
		\begin{assettheorem}
			Riemann conjecture is true			
		\end{assettheorem}
		\begin{proof}
  
  By Theorem Landau [5]
  $$  
	P (s)=ln(\zeta ( s)-Q (s), \,\,  \,Re(s)>1/2+\delta,
	  $$
  let
  $$
  ln\mu(s)=\nu(s) +Q(s)
  $$
  and $ \mu(s)$ another solution (3) and $ \mu(s)$   analitical\,\,\, in $ (s_i < Im(s) <s_{i+1})$ 
  $$
  \nu(s)\big|_{Re (s)>1, s_i < Im(s) <s_{i+1}}=(ln(\zeta(s))-Q(s))\big|_{Re (s)>1, s_i < Im(s) <s_{i+1}}
  $$
    then 
     $$       
  \nu(s)\big|_{Re (s)>1/2, s_i < Im(s) <s_{i+1}}=(ln(\zeta(s))-Q(s))\big|_{Re (s)>1/2, s_i < Im(s) < s_{i+1}}=P (s)
  $$
     
	from last statement  and  	the  estimate Theorem (1-2)  and its symmetry between $\psi_+ $ and $\psi_- $ relative to the critical line leads to the fact that there can be zero zeta-functions
only on the critical line that completes the proof.

		\end{proof}

The results of this paper, together with the results of Voronin's theorem [7] on the universality of the Riemann zeta function, lead to the following applications

\section{Applications}	 
	  \begin{itemize}
   \item  Information compression and information transfer,
   \item  Seismic exploration , 
  \item  Optimal control of oil and gas production  ,
    \item  Optimal control of oil and gas  transfer,
  \item Control of heat and mass transfer in thermonuclear processes,
 \item  Digitalization of technological processes,
\item  Earthquake predictions,
\item Stock market forecast,
\item Stock market forecast,
\item Generations of quantum computers and quantum computing
\item Seismic exploration
\item Tomography
\item Early diagnosis of diseases

\end{itemize}


			{}
        \end{document}